\documentclass[12pt]{amsart}
\usepackage{amsmath,amssymb,amsthm,amsfonts}
\usepackage{setspace,xypic}
\usepackage{graphicx}

\newcommand\inv{^{-1}}

\newcommand\G{\Gamma}

\newcommand\g{\gamma}
\newcommand\Ra{\mathbb R}
\newcommand\Ca{\mathbb C}

\mathchardef\GG="321D

 \DeclareMathOperator{\Isom}{Isom}

\DeclareMathOperator{\Diff}{Diff} 
 
\DeclareMathOperator{\tr}{tr}

 \DeclareMathOperator{\CAT}{CAT}

\newtheorem{theorem}{Theorem}[section]

\newtheorem{lemma}[theorem]{Lemma}

\newtheorem{defn}[theorem]{Definition}

\begin{document}

\title[Superrigidity and harmonic maps]{Cocycle superrigidity and harmonic maps with infinite
dimensional targets}
\author{David Fisher and Theron Hitchman}
\thanks{First author partially supported by NSF grant DMS-0226121 and DMS-0541917.\\
Second author partially supported by NSF VIGRE grant DMS-0240058.}

\begin{abstract}
We announce a generalization of Zimmer's cocycle superrigidity
theorem proven using harmonic map techniques.  This allows us to
generalize many results concerning higher rank lattices to all
lattices in semisimple groups with property $(T)$.  In particular,
our results apply to $SP(1,n)$ and $F_4^{-20}$ and lattices in
those groups.

The main technical step is to prove a very general result concerning
existence of harmonic maps into infinite dimensional spaces, namely
a class of simply connected, homogeneous, aspherical Hilbert
manifolds. This builds on previous work of Corlette-Zimmer and
Korevaar-Schoen. Our result is new because we consider more general
targets than previous authors and make no assumption concerning the
action.  In particular, the target space is not assumed to have
non-positive curvature and in important cases has significant
positive curvature.  We also do not make a ``reductivity" assumption
concerning absence of fixed points at infinity.

The proof of cocycle superrigidity given here, unlike previous
ergodic theoretic ones, is effective.  The straightening of the
cocycle is explicitly a limit of a heat flow.  This explicit
construction should yield further applications.
\end{abstract}

\maketitle

\section{Introduction}
\label{section:intro}

Let $G$ be a semisimple Lie group with all simple factors of real
rank at least two, $\G<G$ a lattice, $S$ a compact manifold and
$\nu$ a volume form on $S$.  In \cite{Z6}, Zimmer proposed a program
to study volume preserving actions of $\G$ on $S$, with the aim of
classifying homomorphisms
$\rho:\G{\rightarrow}\Diff_{\nu}^{\infty}(S)$. A main impetus for
this program was his {\em cocycle superrigidity theorem} which gives
strong information concerning the associated action of $\G$ on any
natural principal bundle over $S$, see \cite{Z1,Z2,Z6,Z7,Z8} and
also \cite{Fe,FM1} for more discussion.

Zimmer's cocycle superrigidity theorem is profoundly influenced in
both statement and proof by Margulis superrigidity theorem for
finite dimensional representations of $\G$ \cite{M1}.  As a
consequence of Margulis' theorem, one has an essentially complete
classification of representations $\G{\rightarrow}GL(n,k)$ where $k$
is a local field. Zimmer's program can be viewed as an attempt to
generalize this problem to consider {\em non-linear representations}
instead.

In the classical setting of Margulis' superrigidity theorems, the
results have been extended by work of Corlette and Gromov-Schoen
\cite{Co2,GS} to cover the case where $G$ is merely a semisimple Lie
group with no compact factors which has property $(T)$ of Kazhdan.
The key addition here is that one is now allowed to have factors of
the type $SP(1,n)$ and $F_4^{-20}$.  To obtain as complete a
classification of linear representations from the work of Corlette
and Gromov-Schoen as obtained by Margulis in \cite{M2}, one also
needs cohomology vanishing theorems due to Matsushima-Murakami and
Raghunathan \cite{MM,Rg} as well as some results from \cite{M2}.
Later, Jost-Yau \cite{JY} and Mok-Siu-Yeung \cite{MSY} showed that
one could reprove Margulis' theorems using harmonic map technology,
at least for cocompact lattices.

Some attempts have been made to use harmonic map techniques to
prove Zimmer's cocycle superrigidity theorem for this broader
class of groups $G$, i.e. for $G$ with no compact factors and
property $(T)$ of Kazhdan, see \cite{CZ,KS3}.  The results in
those papers require additional assumptions on the cocycle which
make them difficult to apply. Some striking applications,
particularly to actions on low dimensional manifolds, are given in
each paper.

In this paper, we prove a complete analogue of Zimmer's theorem for
cocycles satisfying a certain weak integrability condition defined
below.  This allows us to generalize a huge class of results for
actions of higher rank lattices to actions of lattices with property
$(T)$.  Our contribution has two main aspects. First, we allow the
cocycle to have target which is not already assumed semisimple.
Second, we remove conditions concerning the image of the cocycle,
mainly by removing any assumptions concerning fixed points at
infinity. In both places our results are analogous to Zimmer's
result, in higher rank, that the algebraic hull of the cocycle is
reductive \cite{Z8}.  In our setting, the algebraic hull is not a
well-defined invariant, so our proofs do not resemble the ergodic
theoretic proofs very closely, though they are loosely inspired by
the techniques of \cite{M2,Z8}.

The results in \cite{CZ,KS3} can be thought of as extending those of
\cite{Co2,GS} to the infinite dimensional context. The papers
\cite{Co2,GS} apply to representations $\rho:\G{\rightarrow}GL(n,k)$
with the Zariski closure of $\rho(\G)$ reductive and the conditions
needed in \cite{CZ,KS3} are analogous to this.  The arguments in
\cite{MM,Rg,M2} which reduce the case of general linear
representations $\rho$ to the case where the Zariski closure
$\rho(\G)$ is reductive. In the finite dimensional setting, the key
point, which follows from \cite{MM,Rg}, is vanishing of $H^1(\G,V)$
where $V$ is a finite dimensional vector space which is $\G$ module.
While the work we do here is analogous to the reduction in
\cite{MM,Rg,M2}, in our setting the linear problem is replaced by a
non-linear one. It is possible to reduce the problem we study to a
question concerning $H^1(\G,F)$ where $F$ is a topological vector
space which is a non-trivial $\G$ module, but there is no inner
product on $F$ which behaves well under the $\G$ action. This forces
us to use a non-linear method to show that $H^1(\G,F)=0$. In fact,
the formulation in terms of vanishing of cohomology is not
particularly useful in our context. See subsections
\ref{subsection:reductions} and \ref{subsection:vbflow} for more
discussion.

In the remainder of the introduction, we formulate our main results
and some applications.  In subsection \ref{subsection:harmonic}, we
state the result on existence and rigidity of harmonic maps that we
use to prove our cocycle superrigidity theorem.  This is of
independent interest. In subsection \ref{subsection:csr} we describe
our cocycle superrigidity theorem and in subsection
\ref{subsection:apps} we highlight some applications of these
results.  In the second half of this announcement, we sketch the
proofs of our harmonic map results, first describing some reductions
of the problem in subsection \ref{subsection:reductions} and then
outlining the main technical step in subsection
\ref{subsection:vbflow}.

For the remainder of this paper $G$ is a semisimple Lie group with
no compact factors and property $(T)$ of Kazhdan and $\G<G$ is a
lattice.

{\it Acknowledgements:} The first author would like to thank
Mu-tao Wang for several helpful conversation concerning harmonic
maps into infinite dimensional spaces, heat flow and the idea of
zero vanishing theorems.

\subsection{Harmonic map results}
\label{subsection:harmonic}

We let $\tilde M = G/K$ where $K$ is a maximal compact subgroup be
the symmetric space associated to $G$.  Throughout this section $M$
will be a compact locally symmetric space with universal cover
$\tilde M$ and $\pi_1(M)=\G$ will be a cocompact lattice.  The
results stated here probably also hold for $M$ of finite volume, but
there is a fundamental problem of existence of harmonic maps, or
even finite energy maps, which remains unresolved.  See
\cite{JL,JY,Sa} for some partial results in this direction.

We will be considering a harmonic map problem for certain
representation $\rho:\G{\rightarrow}\Isom(X)$ where $X$ is an
aspherical Hilbert manifold. We now describe the class of Hilbert
manifolds we consider, their isometry groups and the condition we
need on $\rho(\G)$.

Let $N$ be a simply connected, aspherical, homogeneous Riemannian
manifold. (Here we reserve the term Riemannian manifold for the
finite dimensional case and refer to the infinite dimensional
analogues as Hilbert manifolds.) We can write $N=H/C$ where
$H=\Isom(N)$ and $C$ is a maximal compact subgroup. We assume that
$C$ does not contain any subgroups which are normal in $H$, and thus
the $H$ action on $H/C$ is effective. To simplify exposition, we
will assume that $H$ is the real points of a linear real algebraic
group, then we can represent $H$ as a subgroup of $GL(m,\Ra)$ for
some large $m$. Let $(S,\mu)$ be a probability measure space. For
$\phi_1,\phi_2:S{\rightarrow}H/C$ measurable, we define

$$d_X(\phi_1,\phi_2)^2=\int_Sd_{N}^2(\phi_1(s),\phi_2(s))d\mu(s)$$

Let $\phi_0:S{\rightarrow}G/K$ be $\phi_0(S)=[C]$ and define
$X:=L^2(S,\mu,H/C)$ to be the set of measurable maps
$\phi:S{\rightarrow}H/C$ with $d(\phi,\phi_0)<\infty$. If $S$ is a
finite set, then $X$ is a finite product $\prod_{s{\in}S}H/C$ where
the metric on each factor is scaled by $\mu(s)$.  In general, we
call $X$ a {\em continuum product} of copies of $H/C$ and $X$ is a
simply connected, complete, homogeneous, aspherical, Hilbert
manifold. If $N=H/C$ has non-positive curvature then so does $X$,
and if $N$ has mixed curvature so does $X$.  We describe a subgroup
of the isometry group of $X$. Let $T:S{\rightarrow}S$ be a measure
preserving transformation and $f:S{\rightarrow}H$ a measurable map
with
\begin{equation}
\label{equation:lnl2} \int_S(\ln^+\|f(s)\|)^2<\infty
\end{equation}

\noindent where here we are taking the matrix norm of $f(s)$ defined
by the embedding of $H$ into $GL(m,\Ra)$. We let
$L^2_{\ln}(S,\mu,H)$ be the set of maps which satisfy equation
$(\ref{equation:lnl2})$. Then we can define an isometry $F$ of $X$
by
$$F(\phi)(s)=f(T{\inv}s)\phi(T{\inv}s)$$
for $\phi{\in}X$.  We denote the group of such isometries by
$\Isom_{\mu}(X)$.  We remark that there is a special case of this
construction where we take $S$ to be a compact Riemannian
manifold, $\mu$ to be the Riemannian volume and
$L^2(S,\mu,SL(n,\Ra)/SO(n))$ to be the measurable, $L^2$ sections
of the $SL(n,\Ra)/SO(n)$ bundle over $S$ associated to the frame
bundle, so $n=\dim(S)$. If $\Gamma$ acts smoothly on $S$
preserving $\mu$ then the derivative of the action yields an
isometric action of $\G$ on $L^2(S,\mu,SL(n,\Ra)/SO(n))$. For this
example, $X$ is usually referred to as the ``space of
$L^2$-Riemannian metrics" on $S$, see for example [KS3]. We prove
the following \cite{FH3}.

\begin{theorem}
\label{theorem:existence} Let $M,\tilde M,\G,X$ as above. Then
given $\rho:\G{\rightarrow}\Isom_{\mu}(X)$, there exists a $\G$
equivariant harmonic map $f:M{\rightarrow}X$.  Furthermore $f$ is
either constant or totally geodesic.
\end{theorem}

\noindent The assumption that $\rho(\G)$ lands in $\Isom_{\mu}(X)$
may at first seem similar to the usual assumption in the theory of
equivariant harmonic maps that the action is {\em reductive}.
However, we are not assuming reductivity in the standard sense. It
is instead a consequence of Theorem \ref{theorem:existence} that the
action must be reductive. In general, $\Isom_{\mu}(X)$ may be a
proper subgroup of the full isometry group of $X$, but there are
important cases where it is the whole isometry group. In particular,
if $H/C$ is an irreducible Riemannian symmetric space of the
non-compact type, then we can show $\Isom(X)=\Isom_{\mu}(X)$. So, in
this case, we are making no assumption at all.

The idea to prove cocycle superrigidity result by studying
harmonic maps into infinite dimensional spaces first appears in
work of Korevaar and Schoen \cite{KS1,KS2,KS3}. The foundations of
a theory of harmonic maps with infinite dimensional target spaces
are laid out in those papers and also, independently and from a
somewhat different viewpoint, in work of Jost
\cite{Jost-Book,Jost-CalcVar,Jost-dirichletforms}.

In \cite{KS2,KS3}, Korevaar and Schoen introduce the notion of a
zero vanishing theorem.  This is a method for proving the
existence of equivariant harmonic maps into non-positively curved
spaces of finite or infinite dimension without assuming
reductivity of the action. Zero vanishing theorems are developed
further in work of Gromov, Izeki-Nayatani and Schoen-Wang
\cite{Gr3,IN,SW}. The general philosophy of zero vanishing
theorems, as well as all the results mentioned, work in a context
where one knows a priori that any equivariant harmonic map must be
constant, and the method in all cases is to produce the resulting
fixed point for the action without first invoking an existence
result. Though our work is similar to these works in that we
produce equivariant harmonic maps without any reductivity
assumption, the method is necessarily quite different because we
need to consider a setting where non-constant harmonic maps exist.
A key step in our proofs, and our primary analytic innovation, is
what one might call a ``relative zero vanishing theorem", see
Theorem \ref{theorem:fiberflow} and the discussion in subsection
\ref{subsection:vbflow}.  This result is the key mechanism for
allowing us to produce non-constant equivariant harmonic maps when
the technique of zero vanishing theorems cannot apply.

There is another point of view, introduced in \cite{Gr2} and used
in \cite{CZ} of studying foliated harmonic maps instead of
harmonic maps with infinite dimensional targets. It is probably
possible to formulate and prove our results from that point of
view as well. For the problems that concern us, the difference in
these points of view is merely whether one thinks of the variable
in $S$ as belonging to the domain or the range. That is, one may
try to solve a family of harmonic map problems parameterized by
$S$ or try to solve a harmonic map problem into a product of
spaces that is parameterized by $S$. Further work using this point
of view to study superrigidity questions was pursued by Benoit
Rivet in an unpublished Ph.D. thesis \cite{Ri}.

\subsection{Cocycle super-rigidity results}
\label{subsection:csr}

In this section we formulate our main result concerning cocycle
super-rigidity.  In what follows, we let $G$ be a semisimple Lie
group with no compact factors and property $(T)$ of Kazhdan,
$\G<G$ a lattice and $(S,\mu)$ a probability measure space.

We recall some basic notions concerning cocycles. Given a group
$D$, a space $S$ and an action $\rho:D{\times}S{\rightarrow}S$, we
define a {\em cocycle over the action} as follows.  Let $L$ be a
group, then a cocycle is a map $\alpha:D{\times}S{\rightarrow}L$
such that $\alpha(g_1{g_2},s)=\alpha(g_1, g_2s)\alpha(g_2,s)$ for
all $g_1,g_2{\in}D$ and all $s{\in}S$. The regularity of the
cocycle is the regularity of the map $\alpha$. If the cocycle is
measurable, we only insist on the equation holding almost
everywhere in $S$. Note that the cocycle equation is exactly what
is necessary to define a skew product action of $D$ on
$S{\times}L$ or more generally an action of $D$ on $S{\times}Y$ by
$d{\cdot}(x,y)=(dx,\alpha(d,x)y)$ where $Y$ is any space with an
$L$ action.

We say two cocycles $\alpha$ and $\beta$ are {\em cohomologous} if
there is a map $\phi:S{\rightarrow}L$ such that
$\alpha(d,s)={\phi(ds)}{\inv}\beta(d,s){\phi(s)}$.  Again we can
define the cohomology relation in any category, depending on how
much regularity we seek or can obtain on $\phi$. A cocycle is
called {\em constant} if it does not depend on $s$, i.e.
$\alpha_{\pi}(d,s)=\pi(d)$ for all $s{\in}S$ and $d{\in}D$.  One
can easily check from the cocycle equation that this forces the
map $\pi$ to be a homomorphism $\pi:D{\rightarrow}L$.  When
$\alpha$ is cohomologous to a constant cocycle $\alpha_{\pi}$ we
will often say that $\alpha$ is cohomologous to the homomorphism
$\pi$.  The cocycle superrigidity theorems imply that many
cocycles are cohomologous to constant cocycles, at least in the
measurable category.

Before stating our results, we require one more definition.

\begin{defn}
\label{defn:L^2} Let $D$ be a locally compact group, $(S,\mu)$ a
standard probability measure space on which $D$ acts preserving
$\mu$  and $L$ a normed topological group. We call a cocycle
$\alpha:D{\times}S{\rightarrow}L$ over the $D$ action {\em $L^2$} if
for any compact subset $M\subset{D}$, the function
$Q_{M,\alpha}(x)=\sup_M\ln^+\|\alpha(m,x)\|$ is in $L^2(S)$.
\end{defn}

Here as above, we take the target of our cocycles, $H$, to be the
real points of a real algebraic group in order to have convenient
norms to work with.  This is not necessary, but does simplify
notations and discussion. We prove the following superrigidity
theorems for cocycles \cite{FH3}.

\begin{theorem}
\label{theorem:Gsuperrigidity} Let $G,S,\mu,H$ be as above. Assume
$G$ acts ergodically on $Y$ preserving $\mu$. Let
$\alpha:G{\times}S{\rightarrow}H$ be an $L^2$, Borel cocycle. Then
$\alpha$ is cohomologous to a cocycle $\beta$ where
$\beta(g,s)=\pi(g){c(g,s)}$. Here ${\pi:G{\rightarrow}H}$ is a
continuous homomorphism and $c:G{\times}S{\rightarrow}C$ is a
cocycle taking values in a compact group centralizing $\pi(G)$.
\end{theorem}

\begin{theorem}
\label{theorem:Gammasuperrigidity} Let $G,\Gamma,S,H$ and $\mu$ be
as above. Assume $\Gamma$ acts ergodically on $S$ preserving $\mu$.
Assume $\alpha:\Gamma{\times}S{\rightarrow}H$ is an $L^2$, Borel
cocycle. Then $\alpha$ is cohomologous to a cocycle $\beta$ where
$\beta(\gamma,s)=\pi(\gamma){c(\gamma,s)}$. Here
${\pi:G{\rightarrow}H}$ is a continuous homomorphism of $G$ and
$c:\Gamma{\times}S{\rightarrow}C$ is a cocycle taking values in a
compact group centralizing $\pi(G)$.
\end{theorem}

The proofs of these results proceed as follows.  For $\G$
cocompact, Theorem \ref{theorem:Gammasuperrigidity} is a direct
consequence of Theorem \ref{theorem:existence}.  We can then
deduce Theorem \ref{theorem:Gsuperrigidity} from Theorem
\ref{theorem:Gammasuperrigidity} using Iozzi's thesis and the
method described in the paper of Corlette and Zimmer \cite{CZ,I}.
It is then standard to prove Theorem
\ref{theorem:Gammasuperrigidity} from Theorem
\ref{theorem:Gsuperrigidity} by inducing actions and cocycles. The
only non-trivial step in that process is proving that the $L^2$
condition is preserved by induction, but this is easy to check in
the cases that concern us.

We remark that the assumption of ergodicity in Theorems
\ref{theorem:Gsuperrigidity} and \ref{theorem:Gammasuperrigidity}
are not exactly necessary, but that without it the statements
become more complicated.  See \cite{FMW} for a method of deducing
the non-ergodic case from the ergodic one, as well as for precise
statements in the non-ergodic case.  In our setting, we prove the
non-ergodic statements directly from Theorem
\ref{theorem:existence}, though the proof of that theorem often
requires that one analyze the actions one ergodic component at a
time.

\subsection{Some applications}
\label{subsection:apps}

In this subsection, we discuss some applications of Theorems
\ref{theorem:Gsuperrigidity} and \ref{theorem:Gammasuperrigidity}.
The gist of this section is that most results concerning smooth
actions of higher rank groups that were proven using super-rigidity
for cocycles now apply more broadly.  This includes all of the
results in \cite{Z1,Z2,Z3,Z4,Z5,Z6,Z7}.

Some notable results, such as the long string of local rigidity
results beginning with \cite{Hu} and culminating in the work
\cite{FM1,FM2,FM3} are not quite established by our methods.  One
step in the proof of these results is not entirely in place: all
methods used for showing smoothness along dynamical foliations
depend on the presence of higher rank abelian subgroups, see
particularly \cite{Hu,KaSp}. While it should be possible to
replace these methods by the methods developed by the second
author in \cite{H}, we are instead pursuing a project to prove
local rigidity directly using the estimates used in the proofs of
Theorems \ref{theorem:existence} and \ref{theorem:h1} and the
criterion for local rigidity introduced in \cite{F}.  We will
discuss this approach elsewhere. We note here that the $C^{3,0}$
local rigidity of \cite[Theorem 1.1]{FM3} is an immediate
consequence of the methods here.

More precisely let $G$ be a (connected) semisimple Lie group with
no compact factors and property $(T)$ of Kazhdan, and $\Gamma<G$
is a lattice. Then combining Theorems \ref{theorem:Gsuperrigidity}
and \ref{theorem:Gammasuperrigidity} with the methods of
\cite{FM1,FM2,FM3} yields:

\begin{theorem}
\label{theorem:main} Let $\rho$ be a quasi-affine action of $G$ or
$\Gamma$ on a compact manifold $X$.   Then the action is $C^{3,0}$
locally rigid.
\end{theorem}

For the reader's benefit, we recall the definition of quasi-affine
from \cite{FM3}.  We remark here that our definition of quasi-affine
implies that all actions considered in Theorem \ref{theorem:main}
are volume preserving.  First recall that, for any Lie group $H$ and
any closed subgroup $\Lambda$, an {\em affine diffeomorphism} $d$ of
$H/{\Lambda}$ is one covered by a diffeomorphism $\tilde d$ of $H$
of the form $\tilde d=A{\circ}T_h$ where $A$ is an automorphism of
$H$ such that $A(\Lambda)=\Lambda$ and $T_h$ is left translation by
$h{\in}H$.

\begin{defn}
\label{definition:affine} {\bf a)} Let $H$ be a connected real
algebraic group, $\Lambda<H$ a cocompact lattice. Assume a
topological group $G$ acts continuously on $H/{\Lambda}$. We say
that the $G$ action on $H/{\Lambda}$ is {\em affine} if every
element of $G$ acts via an affine diffeomorphism.

\noindent {\bf b)} More generally, let $M$ be a compact manifold.
Assume a group  $G$ acts affinely on $H/{\Lambda}$. Choose a
Riemannian metric on $M$ and a cocycle
$\iota:G{\times}H/{\Lambda}{\rightarrow}\Isom(M)$. We call the skew
product action of $G$ on $H/{\Lambda}{\times}M$ defined by
$d{\cdot}(x,m)=(d{\cdot}x, \iota(d,x){\cdot}m)$ a {\em quasi-affine
action}.
\end{defn}

\noindent For $\G$ and $G$ as above, all affine actions are
classified by results in \cite{FM1}.

We note here one other application.  The action on the space of
metrics corresponds to studying the derivative cocycle, and this is
the primary mechanism for many applications.  Here we state results
that use other cocycles over the action to motivate our more general
results.  The result we state generalizes a result due to the first
author and Zimmer.

\begin{theorem}[\cite{FZ}] Let $\G<G$ be a lattice, where $G$ is a noncompact
simple Lie group with property $(T)$ of Kazhdan. Suppose $\G$ acts
analytically and ergodically on a compact manifold $S$ preserving
volume and an analytic connection. Then either:

\begin{enumerate}

\item the action is isometric and $S=K/C$ where $K$ is a compact
Lie group, $C$ is a closed subgroup and the action is by right
translation via $\rho:\G\hookrightarrow K$, a dense image
homomorphism, or

\item there exists an infinite image linear representation
$\sigma:\pi_1(M)\to GL_n(\Ra)$, such that the algebraic
automorphism group of the Zariski closure of
$\sigma\bigl(\pi_1(M)\bigr)$ contains a group locally isomorphic
to $G$.

\end{enumerate}
\end{theorem}

\noindent{\bf Remarks:}\begin{enumerate} \item The analytic
connection in the statement of the theorem can be replaced by any
analytic rigid geometric structure in the sense of \cite{Gr1}. For
$G$ actions a similar result is proved in \cite{Gr1}.  \item Finer
information on the representation in conclusion $(2)$ is obtained
in \cite{FZ}. To obtain the same information here, we need to
extend our results to cover cocycles into groups defined over
other local fields. This is work in progress.
\end{enumerate}

This theorem is proven by studying cocycles associated to linear
representations of the fundamental group of $M$.  Other theorems
proven by studying that type of cocycle are also easily
generalized using our results modulo the remark in $(2)$ above.
See for example \cite{FW,Sch} and references there.

\subsection{Vanishing Cohomology}
\label{section:vanishingcohomology}

Our methods also allow us to deduce new cohomology vanishing results
for certain representations for the groups $\G$ as defined above.
Here $\G$ is a cocompact lattice in $G$ as above. We prove the
following theorem \cite{FH4}.

\begin{theorem}
\label{theorem:h1} Let $\sigma:G\rightarrow GL(V)$ be a finite
dimensional representation of $G$ and
$\pi:\G{\rightarrow}\mathcal{U} (\mathcal{H})$ a unitary
representation of $\G$ on a separable Hilbert space $\mathcal{H}$.
Then $H^1(\G,\sigma{\otimes}\pi)=0$.
\end{theorem}

{\noindent}{\bf Remarks:} \begin{enumerate} \item For certain
representations of $G$ which are of the type described, this result
is known.  Namely it is known for $G$ representations which are
induced from {\em finite dimensional} $\G$ representations of this
kind. These representations are $\G$ representations which induce to
{\em automorphic} representations of $G$. Proofs given in that
context can be made to work somewhat more generally but do not
appear to prove the result stated here, see \cite{BW} and references
there. \item The proof of Theorem \ref{theorem:existence} uses an
estimate that is also used in the proof of Theorem \ref{theorem:h1}.
Theorem \ref{theorem:h1} is proven by using a Bochner type estimate
to give a lower bound on a Laplacian on the associated
$V{\otimes}\mathcal{H}$ bundle over $M$. \item A key step in the
proof of Theorem \ref{theorem:existence} is proving a generalization
of Theorem \ref{theorem:h1} that requires a non-linear heat flow
argument.
\end{enumerate}

\section{On proofs}
\label{section:proofs}

In this section, we discuss the proofs of Theorem
\ref{theorem:existence} and Theorem \ref{theorem:h1}.  In the first
subsection, we describe the relation between the results in
subsection \ref{subsection:harmonic} and those in subsection
\ref{subsection:csr}, as well as stating some results which are used
in the proof. The key point is that cocycles are simply a method for
writing the action on $L^2(S,\mu,H/C)$ in coordinates. In the second
subsection, we sketch the reduction the proof of Theorem
\ref{theorem:existence} to a special case, whose proof we outline in
the final subsection.

For simplicity, throughout this section, we will assume that the
$\G$ action on $(S,\mu)$ is ergodic.  This is tantamount to saying
that there is no non-trivial $\G$ invariant splitting of $X$ as a
product $X=L^2(S_1,\mu|_{S_1},H/C){\times}L^2(S_2,\mu|_{S_2},H/C)$
where $S_1,S_2$ are subsets of $S$ with $S_1{\cup}S_2=S$ up to sets
of measure zero.

\subsection{Cocycles as coordinates for actions on $L^2(S,\mu,H/C)$}
\label{subsection:coordinates}

In this section, we set up the general method of translating between
properties of $H$ valued cocycles over the $\G$ action on $S$ and
properties of the $\G$ action on $X=L^2(S,\mu,H/C)$ coming from a
homomorphism $\rho:D{\rightarrow}\Isom_{\mu}(X)$.

Given $\phi{\in}L^2(S,\mu,H/C)$, our assumption on $\rho$ implies
that each $\rho(\g)$ can be written as

$$(\rho(\g){\phi})(s)=f_{\gamma}(\g{\inv}s)\phi(\g{\inv}s)$$

\noindent where $f:S{\rightarrow}H$ is in $L^2_{\ln}(S,\mu,H)$. The
fact that $\rho$ defines a $\Gamma$ action immediately implies that
the map $\alpha_{\rho}(\g,s)=f_{\gamma}(s)$ is a cocycle over the
$D$ action on $S$.  That $\alpha$ is an $L^2$ cocycle is immediate.
Conversely, given an $L^2$ cocycle $\alpha$, we can define an action
$\rho$ using these equations. Similar statements can be made
concerning action/cocycles of continuous groups, though some care
needs be taken concerning continuity.

The following three lemmas translate geometric properties of
$\rho$ into algebraic and dynamical properties of $\alpha_{\rho}$.
The first is almost trivial and does not really depend on any
assumption on the acting group $D$ or on ergodicity of the $D$
action on $S$.

\begin{lemma}
\label{lemma:fixedpoint} Let
$\rho:D{\rightarrow}\Isom_{\mu}(L^2(S,\mu,H/C))$. Then $\pi$ has a
fixed point if and only if $\alpha_{\rho}$ is $L^2$ cohomologous to
a cocycle taking values in (a conjugate of) $C$.
\end{lemma}

\noindent If the action is not ergodic, the cocycle may take
values in different conjugates of $C$ over different ergodic
components of the action on $S$.

For the remaining two lemmas, the  assumption of ergodicity of the
$D$ action on $X$ is important.  Variants can be stated if the
action is not ergodic, but their formulation is considerably more
complicated.

For the following lemma, we need the assumption that $H$ is a
reductive group.  Variants of this Lemma continue to hold if $H/C$
is replaced by any proper $\CAT(0)$ space.  A variant of this
lemma was discovered independently by Furman and Monod in their
work on cocycle super-rigidity  for products of groups
\cite{FuMo}.

\begin{lemma}
\label{lemma:fpinfinity} Let
$\rho:\G{\rightarrow}\Isom_{\mu}(L^2(S,\mu,H/C))$. There is a $\rho$
invariant equivalence class of geodesic rays if and only if there is
a parabolic subgroup $P<H$ and there is a $\rho$ invariant subspace
of $L^2(S,\mu,H/C)$ isometric to $L^2(S,\mu,P/{C'{\cap}P})$ where
$C'$ is a conjugate of $C$.  This is equivalent to $\alpha_{\rho}$
being $L^2$ cohomologous to a cocycle taking values in $P$.
\end{lemma}

The proof of the final lemma is essentially contained in the paper
of Corlette and Zimmer \cite{CZ}.  For this lemma, we require our
standing assumptions on $\G,G$ and $M$, but do not require that
$H$ is reductive.

\begin{lemma}
\label{lemma:totally geodesic} Let
$\pi:\G{\rightarrow}\Isom_{\mu}(L^2(S,\mu,H/C))$.  There is a
$\rho$ equivariant totally geodesic, harmonic map from $\tilde M$
to $L^2(S,\mu,H/C)$ if and only $\alpha_{\rho}$ is
$L^2$-cohomologous to a cocycle $\beta(\g,s)=\pi(\g){c(\g,s)}$.
Here ${\pi:G{\rightarrow}H}$ is a continuous homomorphism and
$c:\G{\times}S{\rightarrow}C$ is a cocycle taking values in a
compact group centralizing $\pi(G)$.
\end{lemma}

\noindent The proofs of all three lemmas use a standard argument to
play off ergodicity of the $\G$ action on $S$ against tameness of
various actions of algebraic groups.

\subsection{First reduction for the harmonic map
problem} \label{subsection:reductions}

In this subsection we describe how to reduce the proof of Theorem
\ref{theorem:existence} to a special harmonic map problem whose
solution is described in the next subsection. As $H$ is a Lie
group, we know that $H=L{\ltimes}U$ where $U$ is nilpotent and $L$
is reductive. Our assumptions on $H/C$ guarantee that $L$ has no
compact simple factors and that (up to conjugation) $C$ is a
subgroup in $L$. In this context $H/C$ is a fiber bundle with
fiber $U$ over $L/C$ and this makes $X$ into a fiber bundle

$${\xymatrix{{L^2(S,\mu,H/C)}\ar[d]_{\pi}&{F(S,\mu,U)}\ar[l]^{i}\\
{L^2(S,\mu,L/C)}\\}}$$

\noindent where $F(S,\mu,U)=\{f{\in}L^2(S,\mu,H/C)|f(y){\in}U
\text{a.e.}\}$.  We note here that  $L^2(S,\mu,U) \subsetneqq
F(S,\mu,U)$ unless the $L$ action on $U$ is trivial. The method of
proof is to show that the $\G$ action preserves a totally geodesic
copy of $L^2(S,\mu,L/C)$ in  $L^2(S,\mu,H/C)$.  To do this, we can
use the fibered structure of $U$ to inductively reduce to the case
where $U=\Ra^n$. I.e. to studying the case where the above fibration
is replaced by:

$${\xymatrix{{L^2(S,\mu,H/C)}\ar[d]_{\pi}&{F(S,\mu,\Ra^n)}\ar[l]^{i}\\
{L^2(S,\mu,L/C)}\\}}$$

\noindent where
$F(S,\mu,\Ra^n)=\{f{\in}L^2(S,\mu,H/C)|f(y){\in}\Ra^n
\text{a.e.}\}$.   We show that $\rho$ projects to a homomorphism
$\bar \rho:\Gamma{\rightarrow}\Isom(L^2(S,\mu,L/K))$.  We can study
$\rho$ first by using results of Korevaar and Schoen [KS3], namely
the following:

\begin{theorem}[Korevaar-Schoen]
\label{theorem:ks} Let $M,\tilde M, \G$ and $L$ as above and
$X=L^2(S,\mu,L/K)$.  Then either \begin{itemize} \item[\rm{a.}]
there exists a $\G$ equivariant harmonic map $\phi:\tilde
M{\rightarrow}L^2(S,\mu,L/K)$ which is either constant or totally
geodesic, or \item[\rm{b.}] there is an asymptotic equivalence class
$[c]$ of rays in $L^2(S,\mu,L/K)$ invariant under $\G$.
\end{itemize}
\end{theorem}

\noindent{\bf Remarks:}\begin{enumerate} \item Korevaar and Schoen
only claim this theorem in the special case of actions on the
space of metrics.  It is easy to check that their proof goes
through in the generality given here. \item More is true, they
construct a functional $I$, related to energy, which is eventually
decreasing along any ray in the equivalence class $[c]$. \item One
can reformulate the existence of $[c]$ in terms of the existence
of a fixed point on $\partial X$.
\end{enumerate}

Our proof of Theorem \ref{theorem:existence} in the special case
where $H$ is reductive amounts to eliminating case $b.$ of Theorem
\ref{theorem:ks}. To do this, we apply Lemma
\ref{lemma:fpinfinity}, which produces a parabolic subgroup $P<H$
such that $\rho(\G)$ preserves a subspace of the form
$L^2(S,\mu,P/(P{\cap}C'))$ where $C'$ is a conjugate of $C$. We
can, in fact, choose $P$ minimal with this property. This uses the
fact that stabilizers of rays in $L/C$ are exactly (conjugacy
classes of) parabolic subgroups which are ordered by inclusion.
Now $P$ is of the form $L'{\ltimes}U'$ where $L$ is reductive and
$U'$ is unipotent. We can then repeat the discussion above,
replacing $H$ with $P$. However, when we return to apply Theorem
\ref{theorem:ks} again, we are forced to be in case $(a)$ by our
assumption on minimality of $P$.  Some care is required in this
reduction as $L'$ may have compact simple factors and additional
arguments are required in that case.

Combining the reductions we have discussed so far, we can show that
the problem of finding a harmonic map into $L^2(S,\mu,H/C)$ reduces
to the case where this space has the structure:

$${\xymatrix{{L^2(S,\mu,H/C)}\ar[d]_{\pi}&{F(S,\mu,\Ra^n)}\ar[l]^{i}\\
{L^2(S,\mu,L/C)}\\}}$$

\noindent and we assume there is a $\G$ equivariant harmonic map
$f:\tilde M {\rightarrow} L^2(S,\mu,L/C)$. Note that there is an
embedding of $L<H$ which defines an embedding $L/C\subset H/C$ and
an embedding $L^2(S,\mu,L/C){\subset}L^2(S,\mu,H/C)$ and that we can
use this embedding to explicitly trivialize the bundle structure on
$L^2(S,\mu,H/C)$.  Showing that any harmonic map into this bundle
takes values in a translate of $L^2(S,\mu,L/C)$ is formally
equivalent to vanishing of $H^1(\G,F(S,\mu,\Ra^n))$ where
$F(S,\mu,\Ra^n)$ becomes a $\Gamma$ module by pulling back the
vector bundle over $\tilde M$ by $f$ and then noting that it
descends to a vector bundle on $M$ by $\G$ equivariance of $f$. The
$\G$ module structure on $F(S,\mu,\Ra^n)$ can be given more
explicitly, since under our assumptions, one can write the $\G$
action on $L^2(S,\mu,H/C)$ as
$$\g{\cdot}(\phi_1(s),\phi_2(s))=(\beta(\g,s)\phi_1(s),\beta(\g,s)\phi_2(s)+h(\g,s))$$
\noindent where $\phi=(\phi_1,\phi_2)$ with $\phi_1$ taking values
in $L$ and $\phi_2$ takes values in $\Ra^n$. We know that
$\beta(\g,s)=\sigma(\g)c(\g,s)$ where $\sigma$ is a $G$
representation.  Note that the representation $\tau$ of $\G$ on
$L^2(S,\mu,\Ra^n)$ defined by
$(\g{\cdot}{\phi})(s)=c(\g,\g{\inv}s)\phi(\g{\inv}s)$ is unitary.
Here $h$ is an $F(S,\mu,\Ra^n)$ valued cocycle over the linear
action of $\G$ determined by $\beta$.

If $F(S,\mu,\Ra^n)=L^2(S,\mu,\Ra^n)$ or if $h{\in}L^2(S,\mu,\Ra^n)$
then what is needed is exactly Theorem \ref{theorem:h1}. However, it
is only the case that $F(S,\mu,\Ra^n)=L^2(S,\mu,\Ra^n)$ when the map
$f$ is constant, in which case, vanishing of
$H^1(\G,L^2(S,\mu,\Ra^n))$ is an immediate consequence of property
$(T)$ since the $\G$ representation on $L^2(S,\mu,\Ra^n)$ is unitary
in this case. When $f$ is not constant, the fact that the distance
we have on $F(S,\mu,\Ra^n)$ is not induced by an inner product
structure but by the embedding of $F(S,\mu,\Ra^n)$ in $
L^2(S,\mu,\Ra^n)$ introduces profound difficulties. As a result we
need to use a heat flow applied to the maps into $L^2(S,\mu,H/C)$
which project to $f$.  We describe this in the next section.

\subsection{Heat flow for maps into vector bundles}
\label{subsection:vbflow}

In this section, we describe a heat flow method which we apply to
complete the proof of Theorem \ref{theorem:existence}.   Here as
above, we assume that $\G,G$ and $M$ are as in Theorem
\ref{theorem:existence}.  For the reader's convenience, we
explicitly state the result we prove here:

\begin{theorem}
\label{theorem:fiberflow} Let $H=L{\ltimes}\Ra^n$ be a Lie group
where $L$ is reductive and has no compact factors and $C<L$ is a
maximal compact subgroup. Let
$\rho:\G{\rightarrow}\Isom_{\mu}(L^2(S,\mu,H/C))$ be a homomorphism
which projects to $\bar
\rho:\G{\rightarrow}\Isom_{\mu}(L^2(S,\mu,L/C))$ and let $f:\tilde
M{\rightarrow}L^2(S,\mu,L/C)$ be a $\bar \rho$-equivariant, totally
geodesic, harmonic map. Then there is a lift of $f$ to a map $\tilde
f:M{\rightarrow}L^2(S,\mu,H/C)$ such that $\tilde f$ is harmonic,
$E(\tilde f)=E(f)$ and $\tilde f$ takes values in some translate of
$L^2(S,\mu,L/C)$ in $L^2(S,\mu,H/C)$.
\end{theorem}

\noindent{\bf Remarks:}\begin{enumerate} \item We assume $f$ is
totally geodesic, since the case of $f$ constant is easy by the
final remarks of the last subsection. \item It follows from Theorem
\ref{theorem:fiberflow} and a standard argument using Bochner
estimates of Corlette \cite{Co2}, Jost-Yau \cite{JY} or
Mok-Siu-Yeung \cite{MSY} that $\tilde f$ is totally geodesic.
\item This theorem can be reformulated to say that
$H^1(\G,F(S,\mu,\Ra^n))=0$, but this point of view is not
particularly helpful for understanding the proof.
\end{enumerate}

We study this harmonic map problem by analyzing a heat flow that we
construct here.  The space on which the heat flow is constructed is
$$\mathcal F=\{g:G/K \rightarrow L^2(S,\mu,H/C) |\ g \text{ is } \rho \text{-equivariant and }
\pi{\circ}g=f\}.$$ \noindent If we let $\tau(g)=-\tr\nabla dg$ be
the tension field of $g$ with respect to the natural metric on
$TM^*{\otimes}TY$, then we can compute
$\tau(g)=(\tau_1(g),\tau_2(g))$ where our coordinates are from the
splitting $g^*TY=g^*TY'{\oplus}g^*TF$. It is straightforward to
check that the condition $\pi{\circ}g=f$ implies that $\tau_1(g)=0$.
As noted above, $X$ admits a global trivialization as
$Y=L^2(S,\mu,L/C){\times}F(S,\mu,\Ra^n)$ in which we can write
$g(m)=(f(m),g_F(m))$. Therefore the tension field is vertical with
respect to our fibration. This is what allows us to construct a heat
flow that preserves the class $\mathcal F$ and which only alters
$g_F$.

The key point for the rest of the argument is to understand the
structure of $TM^*{\otimes}g_F^*TF$.  First $TF$ is a trivial
$L^2(S,\mu,\Ra^n)$ bundle over $F$ and the associated bundle
$g_{F}^*L^2(S,\mu,\Ra^n)$ is easily seen to be a flat bundle where
the holonomy is a representation of the form $\sigma{\otimes}\tau$
where $\sigma$ is the restriction of a finite dimensional $G$
representation and $\tau$ is a unitary $\G$ representation. In fact,
as described in the last subsection we can  write the action of $\G$
on $X$ in coordinates and then the $\G$ module structure on the
bundle $g_F^*L^2(S,\mu,\Ra^n)$ is given explicitly.

Short time existence of heat flow follows from an adaptation of
standard arguments using implicit function theorems. Long time
existence is then shown using the fact that $g^*L^2(X,\mu,\Ra^n)$ is
a flat bundle and the Eells-Sampson Bochner formula in a manner
similar to that in \cite{Co1}.  Both of these steps are, of course,
complicated by the fact that our target space is not finite
dimensional.

To understand the long time behavior of the heat flow requires use
of finer structure.  In particular, we use another Bochner type
estimate to show that the total energy of $g_2$ goes to zero as
$t{\rightarrow}\infty$.  Once we know this, Theorem
\ref{theorem:fiberflow} is immediate, since this implies that
$g_2$ is constant.

The Bochner-type estimate we use in this context is the same one
that is used to prove Theorem \ref{theorem:h1}.  In the context of
Theorem \ref{theorem:h1}, this estimate is used to produce a lower
bound on the Laplacian on one forms of the form $(\Delta
u,u)>c\|u\|$ for some absolute constant $c$.  This forces cohomology
to vanish by \cite[Proposition 1.3.1]{Mk}.  The same estimate
applied in the non-linear setting implies $\|\tau(g_2)\|>cE(g_2)$.
If $g_2^t$ is the heat flow applied to $g_2$ at time $t$, then this
implies that
$$\frac{dE(g_2^t)}{dt}=-\|\tau(g^t_2)\|<-cE(g^t_2)$$
\noindent which suffices to show that $E(g_2^t)$ decreases
exponentially quickly to zero.

The estimate we need follows by comparing two Laplacians on the
$L^2(X,\mu,\Ra^n)$ bundle on $M$ or, more precisely, by taking the
difference of Bochner formulae for these two Laplacians. Each
Laplacian is associated to a connection. The first connection
$\nabla_1$ is simply the standard flat bundle connection with
associated differential $d_1$ and Laplacian $\Delta_1$. To discuss
the second Laplacian we make a simplifying assumption. Namely, we
assume that $\Ra^n$ is an irreducible $G$ module.  In the full
proof, algebraic arguments are used to reduce to this case. In
this setting, we have a flat $\Ra^n$ bundle over $M$ with holonomy
$\sigma|_{\Gamma}$ and we can build an isomorphic $\Ra^n$ bundle
over $M$ which is associated to the $K$ bundle $G/{\Gamma}$ over
$M$ as in \cite{MM}.  We show that this same process can be used
to give a different structure on the $L^2(S,\mu,\Ra^n)$ bundle.
This other structure on the bundle gives rise to a second
connection $\nabla_2$ with associated differential $d_2$ and
Laplacian $\Delta_2$.  Note that each Laplacian depends on a
choice of metric on $L^2(S,\mu,\Ra^n)$ which depends on the choice
of metric on $\Ra^n$. To have the required estimate, we need the
metric on $\Ra^n$ to be {\em adapted} in the sense of Matsushima
and Murakami \cite{MM}.  If $G_{\Ca}$ is the complexification of
$G$ and $\sigma_{\Ca}$ is the resulting representation on $\Ca^n$,
then a metric on $\Ra^n$ is adapted if it is the restriction of a
Hermitian metric on $\Ca^n$ which is invariant under the compact
form $G_c$ of $G$ in $G_{\Ca}$. The metric on $\Ra^n$ induced by
an invariant metric on $H/C$ is easily checked to have this
property and it is also easy to check that this property is
preserved by pulling back by $g_F$.  If one takes the general
Eells-Sampson Bochner formulae for $\Delta_1$ and $\Delta_2$ on
$1$-forms then a computation similar to the one in \cite{MM}
implies that the difference $\Delta_1-\Delta_2$ is bounded below
by a certain algebraic operator.  This operator is shown to be
strictly positive by Raghunathan in \cite{Rg}, which then forces
$\Delta_1$ to be strictly positive and completes our proofs.  The
interpretation of the Matsushima-Murakami Bochner formula as a
difference of Eells-Sampson type Bochner formulas is implicit in
\cite{MM}, but does not appear to have been noted explicitly
before now.

By taking more care in the reductions than we have done here, it
is possible to see that our argument shows that the functional $I$
on the space $X$ as defined by Korevaar and Schoen attains a
minimum in $X$.  Using this we can show that the harmonic map
obtained here is a limit of a heat flow.

\noindent Department of Mathematics, Indiana University, Rawles
Hall, Bloomington, IN, 47405.\

\noindent Department of Mathematics, Rice University, Houston, TX,
77005.\

\end{document}